ANNALES
DE L'INSTITUT
HENRI
POINCARÉ
PROBABILITÉS
ET STATISTIQUES

www.imstat.org/aihp

# Comparison between criteria leading to the weak invariance principle

## Olivier Durieu and Dalibor Volný

*Laboratoire de Mathématiques Raphaël Salem, UMR 6085 CNRS–Université de Rouen, Avenue de l'Université, BP.12,*
*F76801 Saint-Étienne-du-Rouvray, France. E-mail: olivier.durieu@etu.univ-rouen.fr; dalibor.volny@univ-rouen.fr*

**Abstract.** The aim of this paper is to compare various criteria leading to the central limit theorem and the weak invariance principle. These criteria are the martingale-coboundary decomposition developed by Gordin in *Dokl. Akad. Nauk SSSR* **188** (1969), the projective criterion introduced by Dedecker in *Probab. Theory Related Fields* **110** (1998), which was subsequently improved by Dedecker and Rio in *Ann. Inst. H. Poincaré Probab. Statist.* **36** (2000) and the condition introduced by Maxwell and Woodroofe in *Ann. Probab.* **28** (2000) later improved upon by Peligrad and Utev in *Ann. Probab.* **33** (2005). We prove that in every ergodic dynamical system with positive entropy, if we consider two of these criteria, we can find a function in $\mathbb{L}^2$ satisfying the first but not the second.

**Résumé.** Le but de cet article est de comparer différents critères conduisant au théoreme limite centrale et au principe d'invariance faible. Ces critères sont la décomposition martingale-cobord développée par Gordin dans *Dokl. Akad. Nauk SSSR* **188** (1969), le critère projectif introduit par Dedecker dans *Probab. Theory Related Fields* **110** (1998), par la suite amélioré par Dedecker et Rio dans *Ann. Inst. H. Poincaré Probab. Statist.* **36** (2000) et la condition introduite par Maxwell et Woodroofe dans *Ann. Probab.* **28** (2000), plus tard améliorée par Peligrad et Utev dans *Ann. Probab.* **33** (2005). On montre que dans tout système dynamique ergodique strictement positive, si l'on considère deux de ces critères, on peut trouver une fonction dans $\mathbb{L}^2$ vérifiant le premier mais pas le deuxième.

*MSC:* 60F05; 60F17; 60G10; 28D05; 60G42

*Keywords:* Stationary process; Central limit theorem; Weak invariance principle; Martingale approximation; Projective criterion

## 1. Introduction

Let $(\Omega, \mathcal{A}, \mu)$ be a probability space and $T : \Omega \to \Omega$ a bijective bimeasurable transformation preserving the measure $\mu$ (i.e., $\mu(T^{-1}A) = \mu(A)$, $\forall A \in \mathcal{A}$). $(\Omega, \mathcal{A}, \mu, T)$ is called a dynamical system. We will assume that it is ergodic, i.e., $T^{-1}A = A$ implies $\mu(A) = 0$ or 1. Let $f$ be a measurable function defined on $\Omega$, then $(f \circ T^i)_{i \in \mathbb{Z}}$ is a stationary process. On the other hand, for every stationary random process $(X_i)_{i \in \mathbb{Z}}$, there exists a dynamical system $(\Omega, \mathcal{A}, \mu, T)$ and a function $f$ on $\Omega$ such that $(f \circ T^i)_{i \in \mathbb{Z}}$ and $(X_i)_{i \in \mathbb{Z}}$ have the same distribution (see, e.g., [4], p. 178). We assume that $\mathbb{E}(f) = 0$.

Let $S_n(f) = \sum_{i=0}^{n-1} f \circ T^i$. We say that $f$ satisfies the Central Limit Theorem (CLT) if $\frac{1}{\sqrt{n}} S_n(f)$ converge in distribution to a normal law.

Let $S_n(f,t) = S_{\lfloor tn \rfloor}(f) + (tn - \lfloor tn \rfloor) f \circ T^{\lfloor tn \rfloor}$, where $\lfloor x \rfloor$ denotes the greatest integer that is smaller than $x$. We say that $f$ satisfies the weak invariance principle (or Donsker invariance principle) if the process $\{\frac{1}{\sqrt{n}} S_n(f,t) \mid t \in [0,1]\}$ converges in distribution to a Brownian motion in the space $C[0,1]$ with the uniform norm. In the sequel, we shall call this the invariance principle.

These two limit theorems have been extensively studied and several methods of proving them have been developed. In this paper, we restrict our attention to three of them.

*Martingale-coboundary decomposition*

This method of proving the CLT was first used by Gordin [8]. The idea is to represent $f$ in the form

$$f = m + g - g \circ T,$$

where $(m \circ T^i)_{i \in \mathbb{Z}}$ is a martingale difference sequence. The term $g - g \circ T$ is called a coboundary and $g$ is the transfer function. This decomposition is called a martingale-coboundary decomposition.

If $m \in \mathbb{L}^2(\Omega)$, the CLT for martingale differences of Billingsley [1] and Ibragimov [12] applies. If $g$ is measurable, the telescopic sum $\frac{1}{\sqrt{n}} \sum_{i=0}^{n-1} (g - g \circ T) \circ T^i$ goes to zero in probability. So, if we can find the above decomposition with $m \in \mathbb{L}^2(\Omega)$ and $g$ measurable, the CLT holds for $f$ by application of Theorem 4.1 of [2]. Moreover, in this case, if $g \in \mathbb{L}^2(\Omega)$, we also have the invariance principle, as proved in [11] (see also [10]). On the other hand, there exist counterexamples with $g \in \mathbb{L}^1(\Omega)$ and $g - g \circ T \in \mathbb{L}^2$ where the invariance principle does not hold, see [21]. According to [10, 21], if $m \in \mathbb{L}^2(\Omega)$, a necessary and sufficient condition to have the invariance principle is

$$\frac{1}{\sqrt{n}} \max_{i \le n} |g \circ T^i| \underset{n \to \infty}{\longrightarrow} 0 \quad \text{in probability.}$$

We say that $(f \circ T^i)_{i \in \mathbb{Z}}$ (or $f$) admits a martingale-coboundary decomposition in $\mathbb{L}^p$, $p \ge 1$, if $m$ and $g$ are in $\mathbb{L}^p$. Let $\mathcal{F} \subset \mathcal{A}$ be a $T$-invariant $\sigma$-algebra, i.e., $\mathcal{F} \subset T^{-1}\mathcal{F}$. Note $\mathcal{F}_i = T^{-i}\mathcal{F}$. If we assume that $f$ is $\mathcal{F}_\infty$-measurable and $\mathbb{E}(f|\mathcal{F}_{-\infty}) = 0$, then $f$ admits a martingale-coboundary decomposition in $\mathbb{L}^p$ with $(m \circ T^i)_{i \in \mathbb{Z}}$ adapted to the filtration $(\mathcal{F}_i)_{i \in \mathbb{Z}}$ if and only if the series

$$\sum_{i=0}^{\infty} \mathbb{E}(f \circ T^i | \mathcal{F}_0) \quad \text{and} \quad \sum_{i=0}^{\infty} (f \circ T^{-i} - \mathbb{E}(f \circ T^{-i} | \mathcal{F}_0))$$

converge in $\mathbb{L}^p$, see [10, 20]. This is the characterisation that we shall always use.

Moreover, when the filtration $(\mathcal{F}_i)_{i \in \mathbb{Z}}$ is adapted to the process $(f \circ T^i)_{i \in \mathbb{Z}}$, the second sum equals zero. This will be the case in the sequel.

According to what precedes, existence of martingale-coboundary decomposition in $\mathbb{L}^p$, with $p \ge 2$ implies the invariance principle.

This method gives results in various situations. An interesting example is its application to differentiable dynamical systems. It is well adapted to the hyperbolic case (e.g. [14]), or the partially hyperbolic case (e.g. [13]).

*Projective criterion*

Another method is to establish a projective property developed by Dedecker [5]. He introduced this criterion to prove CLT for random vector fields. Dedecker and Rio [6] have shown that it gives a powerful criterion for proving the invariance principle (in dimension one). We say that $(f \circ T^i)_{i \in \mathbb{Z}}$ (or $f$) satisfies the projective criterion if

$$\sum_{k=1}^{\infty} f \mathbb{E}(f \circ T^k | \mathcal{F}_0) \quad \text{converges in } \mathbb{L}^1,$$

where $(\mathcal{F}_i)_{i \in \mathbb{Z}}$ is a filtration adapted to $(f \circ T^i)_{i \in \mathbb{Z}}$.

According to [6], if $f \in \mathbb{L}^2$ satisfies the projective criterion, then $f$ satisfies the invariance principle (and the CLT).

*Maxwell–Woodroofe condition*

We say that $(f \circ T^i)_{i \in \mathbb{Z}}$ (or $f$) satisfies the Maxwell–Woodroofe condition if

$$\sum_{n=1}^{\infty} \frac{\|\mathbb{E}(S_n(f)|\mathcal{F}_0)\|_2}{n^{3/2}} < \infty.$$

This criterion was first introduced by Maxwell and Woodroofe [15]. They proved that the CLT holds under this condition. Recently, Peligrad and Utev [18] have shown that the same condition also implies the invariance principle.

For examples of applications of the last two methods, the reader can see [16].

Our purpose is to compare the dependence between these criteria. Section 2 contains the statement of our main result while the remainder of the paper is devoted to its proof. Sections 3 and 4 present a general type of a suitable function in a dynamical system. In Section 5, this model is used to produce specific functions proving our result.

## 2. Main results

It is of interest to know whether one of the considered criteria implies another. This is the question that we propose to answer. First, note that a simple application of the Hölder inequality

$$\left\| \sum f\mathbb{E}(f \circ T^k|\mathcal{F}_0) \right\|_1 \leq \left\| \sum \mathbb{E}(f \circ T^k|\mathcal{F}_0) \right\|_p \|f\|_q \tag{1}$$

with $\frac{1}{p} + \frac{1}{q} = 1$, leads to the following remark.

**Remark 1.** *The martingale-coboundary decomposition in $\mathbb{L}^2$ implies the projective criterion.*

**Proof.** It follows from the convergence of $(\sum_{k=1}^{n} \mathbb{E}(f \circ T^k|\mathcal{F}_0))_{n \geq 1}$ in $\mathbb{L}^2$ and the inequality (1) with $p = q = 2$.  □

**Remark 2.** *The martingale-coboundary decomposition in $\mathbb{L}^2$ implies the Maxwell–Woodroofe condition.*

So we are interested in the martingale-coboundary decomposition in $\mathbb{L}^1$. The same kind of arguments show:

**Remark 3.**

(a) *For bounded functions, the martingale-coboundary decomposition in $\mathbb{L}^1$ implies the projective criterion.*
(b) *For a function $f$ such that $|f| > C > 0$, the inverse implication is true.*

**Proof.** (a) follows from application of (1) with $p = 1$ and $q = \infty$. For (b), it is enough to note that $\frac{1}{f}$ is a bounded function.  □

We will see, by counterexamples in $\mathbb{L}^2$, that in general, the martingale-coboundary decomposition in $\mathbb{L}^1$, the projective criterion, and the Maxwell–Woodroofe condition do not result from each other, even if the function verifies the CLT or the invariance principle. Clearly, for the example constructed in [21], which verifies the martingale-coboundary decomposition in $\mathbb{L}^1$ but not the invariance principle, the projective criterion and the Maxwell–Woodroofe condition do not hold. We can also find counterexamples in the class of functions satisfying the invariance principle. Our main result is the following theorem.

**Theorem.** *Let* $(\Omega, \mathcal{A}, \mu, T)$ *be an ergodic dynamical system with positive entropy. In each case, there exists a function in* $\mathbb{L}^2(\Omega)$ *satisfying:*

(i) *the projective criterion but not the martingale-coboundary decomposition in* $\mathbb{L}^1$;
(ii) *the martingale-coboundary decomposition in* $\mathbb{L}^1$ *and the invariance principle, but not the projective criterion;*
(iii) *the Maxwell–Woodroofe condition but not the martingale-coboundary decomposition in* $\mathbb{L}^1$;
(iv) *the martingale-coboundary decomposition in* $\mathbb{L}^1$ *and the invariance principle, but not the Maxwell–Woodroofe condition;*
(v) *the Maxwell–Woodroofe condition but not the projective criterion;*
(vi) *the projective criterion but not the Maxwell–Woodroofe condition.*

## 3. Preliminary

To prove the theorem, in each case, we will produce a function in $\mathbb{L}^2$ satisfying the first condition but not the second one. These functions will be defined in the same way, so we begin by a general construction. The first step is to choose disjoint sets having a nice property. This section is devoted to the exposition of the construction of these sets.

Let $\mathcal{C}$ be a sub-$\sigma$-algebra of $\mathcal{A}$ such that $T^{-1}\mathcal{C} = \mathcal{C}$. We assume that the measure $\mu$ restricted to $\mathcal{C}$ is non-atomic. The goal is to establish Lemma 2 corresponding to the construction of disjoint sets $A_k$ quasi-invariant under a finite number of iterations of the transformation. Moreover, we want to control the measure of the $A_k$. First, we recall the following lemma. A proof can be found in [7], as a particular case of Theorem 2.2. It can also be done directly by using the Rokhlin lemma.

**Lemma 1.** *Let* $N \in \mathbb{N}$, $0 < \rho < 1$ *and* $\varepsilon > 0$. *There exists a set* $A \in \mathcal{C}$ *such that* $\mu(A) = \rho$ *and for all* $i, j \in \{0, \ldots, N\}$,

$$\mu(T^{-i}A \Delta T^{-j}A) \leq \varepsilon.$$

**Remark 4.** *If* $\mu(A_1 \Delta A_2) \leq \varepsilon_1$ *and* $\mu(B_1 \Delta B_2) \leq \varepsilon_2$, *then* $\mu((A_1 \setminus B_1)\Delta(A_2 \setminus B_2)) \leq \varepsilon_1 + \varepsilon_2$.

We are going to use this remark as well as Lemma 1 to show:

**Lemma 2.** *Let* $(N_k)_{k \in \mathbb{N}} \subset \mathbb{N}$ *with* $N_k \nearrow \infty$.
 *Let* $\rho_k = \lambda^k$ $(0 < \lambda < \frac{1}{2})$ *and* $a = 1 - \sum_{k \geq 1} \rho_k \in (0, 1)$.
 *Let* $(\varepsilon_k)_{k \in \mathbb{N}}$ *be a strictly decreasing sequence of positive reals converging to zero.*
 *There exists* $(A_k)_{k \in \mathbb{N}} \subset \mathcal{C}$ *such that:*

(i) *the sets* $A_k$ *are mutually disjoint;*
(ii) $a\rho_k \leq \mu(A_k) \leq \rho_k$, *for all* $k$;
(iii) *for all* $k \geq 1$ *and for all* $i, j \in \{0, \ldots, N_k\}$, $\mu(T^{-i}A_k \Delta T^{-j}A_k) \leq \varepsilon_k$.

**Proof.** First, denote $\delta_k = \varepsilon_k - \varepsilon_{k-1}$, $k \geq 1$.
For every $k \geq 1$, by Lemma 1, there exists a set $A'_k \in \mathcal{C}$ such that $\mu(A'_k) = \rho_k$ and $\mu(T^{-i}A'_k \Delta T^{-j}A'_k) \leq \delta_k$ for all $i, j \in \{0, \ldots, N_k\}$.
We define $A_k = A'_k \setminus \bigcup_{j=k+1}^{\infty} A'_j$.
Hence, (i) holds by construction. We have $\mu(A_k) \leq \mu(A'_k) = \rho_k$ and

$$\mu(A_k) \geq \mu(A'_k) - \sum_{j=k+1}^{\infty} \mu(A'_j) \geq \rho_k - \rho_k \sum_{j=1}^{\infty} \rho_j = a\rho_k.$$

So, (ii) is verified. For (iii), we use the preceding remark to have, for all $k \geq 1$, $i, j \in \{0, \ldots, N_k\}$,

$$\mu(T^{-i}A_k \Delta T^{-j}A_k) \leq \sum_{j=k}^{\infty} \delta_j = \varepsilon_k.$$

$\square$

An important feature of Lemma 2 is that there is no dependence between $N_k$ and $\rho_k$ and the $\varepsilon_k$ can be chosen arbitrarily small.

## 4. General approach

Here, we give the general model from which the counterexamples will be constructed, proving our theorem. We will define a "pattern function" depending on sequences $(N_k)_{k \in \mathbb{N}} \subset \mathbb{N}$, $(\theta_k)_{k \in \mathbb{N}} \subset \mathbb{R}_+$, $(\rho_k)_{k \in \mathbb{N}} \subset (0, 1)$ and $(\varepsilon_k)_{k \in \mathbb{N}}$ ($\varepsilon_k \ll \rho_k$). In Section 5, we will see that changing the values of the sequences provides different counterexamples.

### 4.1. The model

$(\Omega, \mathcal{A}, \mu, T)$ is an ergodic dynamical system with positive entropy. By the Sinai theorem, it admits a factor which is a Bernoulli shift with the same entropy (see [19]). So, it is sufficient to consider the case where $(\Omega, \mathcal{A}, \mu, T)$ is a Bernoulli shift with positive entropy. This means that:

$\Omega = \{0, 1, \ldots, l\}^{\mathbb{Z}}$, for some $l \in \mathbb{N}^* = \{1, 2, \ldots\}$;

$\mathcal{A}$ is the product $\sigma$-algebra;

$\mu$ is the product measure given by $\mu(\{\omega \in \Omega : w_0 = i\}) = p_i$, for $i = 0, \ldots, l$, with $p_i > 0$ and $\sum_{i=0}^{l} p_i = 1$;

$T$ is the left shift on $\Omega$, i.e., $(Tw)_i = w_{i+1}$.

Now, using the Ornstein isomorphism theorem (see [17]), we can see that a Bernoulli shift is isomorphic to a product of two Bernoulli shifts. In particular, our system admits two independent Bernoulli factors. We denote by $\mathcal{B}$ and $\mathcal{C}$ the $T$-invariant $\sigma$-algebras corresponding to them. In order to simplify some proof, we assume that the first one is a Bernoulli $(\frac{1}{2}, \frac{1}{2})$. The reader can check that all the upcoming proofs remain valid for another Bernoulli shift. So, we can define a $\mathcal{B}$-measurable function $e_0 : \Omega \to \{-1, 1\}$ such that $\mu(\{e_0 = -1\}) = \mu(\{e_0 = 1\}) = \frac{1}{2}$ and if $e_i = e_0 \circ T^i$ for $i \in \mathbb{Z}$, then $(e_i)_{i \in \mathbb{Z}}$ is an i.i.d. sequence. Of course, $(e_i)_{i \in \mathbb{Z}}$ is independent of $\mathcal{C}$.

Let $\mathcal{F}_0 = \mathcal{C} \vee \sigma\{e_i \mid i \leq 0\}$ and $\mathcal{F}_k = T^{-k}\mathcal{F}_0 = \mathcal{C} \vee \sigma\{e_i \mid i \leq k\}$, $k \in \mathbb{Z}$.

By application of Lemma 2, we consider the sets $A_k \in \mathcal{C}$ corresponding to sequences $(N_k)_{k \in \mathbb{N}}$, $(\rho_k)_{k \in \mathbb{N}}$ and $(\varepsilon_k)_{k \in \mathbb{N}}$. The function $f$ is defined by

$$f = \sum_{k=1}^{\infty} f_k \mathbb{1}_{A_k} \quad \text{with } f_k = \theta_k e_{-N_k}, \tag{2}$$

where $\mathbb{1}_A$ is the indicator function of $A$. The $\varepsilon_k$ can be chosen arbitrarily small. So, we shall not define them in each example. We just assume that

$$\sum_{k=1}^{\infty} \theta_k N_k \sqrt{\varepsilon_k} < \infty, \tag{3}$$

which implies $\sum_{k=1}^{\infty} \theta_k N_k \varepsilon_k < \infty$ ($\varepsilon_k < 1$).

We consider the stationary process $(f \circ T^i)_{i \in \mathbb{Z}}$ for which $(\mathcal{F}_i)_{i \in \mathbb{Z}}$ is an adapted filtration.

**Proposition 1.** *The function $f$ belongs to $\mathbb{L}^2$ if and only if $\sum_{k=1}^{\infty} \theta_k^2 \rho_k < \infty$.*

**Proof.** By disjointness of the sets $A_k$,

$$\|f\|_2^2 = \sum_{k=1}^{\infty} \|f_k \mathbb{1}_{A_k}\|_2^2 = \sum_{k=1}^{\infty} \theta_k^2 \mu(A_k).$$

Now, by Lemma 2, $a\rho_k \leq \mu(A_k) \leq \rho_k$. Thus,

$$a \sum_{k=1}^{\infty} \theta_k^2 \rho_k \leq \|f\|_2^2 \leq \sum_{k=1}^{\infty} \theta_k^2 \rho_k. \qquad \square$$

In what follows, we apply the three studied criteria to our function $f$. We express $f$ satisfying one of them by conditions concerning the sequences $(N_k)_{k \in \mathbb{N}}$, $(\theta_k)_{k \in \mathbb{N}}$ and $(\rho_k)_{k \in \mathbb{N}}$.

**Proposition 2.** *The stationary process $(f \circ T^i)_{i \in \mathbb{Z}}$ admits a martingale-coboundary decomposition in $\mathbb{L}^1(\Omega)$ if and only if $\sum_{k=1}^{\infty} \theta_k \sqrt{N_k} \rho_k < \infty$.*

**Proof.** Recall that the function $f$ admits a martingale-coboundary decomposition in $\mathbb{L}^1$ if and only if $(\sum_{i=1}^{n} \mathbb{E}(f \circ T^i | \mathcal{F}_0))_{n \geq 1}$ converges in $\mathbb{L}^1$.

**Necessary condition.** We assume that $\sum_{k=1}^{\infty} \theta_k \sqrt{N_k} \rho_k = \infty$. We shall show that if $\sum_{i=1}^{\infty} \mathbb{E}(f \circ T^i | \mathcal{F}_0)$ converges in $\mathbb{L}^1$ then $\mathbb{E}|\sum_{i=1}^{\infty} \mathbb{E}(f \circ T^i | \mathcal{F}_0)| = \infty$, a contradiction.

For all $k$ and $i$, $\mathbb{1}_{A_k} \circ T^i$ is $\mathcal{F}_0$-measurable, so

$$\sum_{i=1}^{\infty} \mathbb{E}(f \circ T^i | \mathcal{F}_0) = \sum_{i=1}^{\infty} \sum_{k=1}^{\infty} \mathbb{E}(f_k \circ T^i | \mathcal{F}_0) \mathbb{1}_{A_k} \circ T^i.$$

We will use the fact that the measure of $A_k \Delta T^{-i} A_k$ is small when $i \leq N_k$ to simplify the summation. Note that $\mathbb{E}(e_i | \mathcal{F}_0) = e_i$ if $i \leq 0$ and $\mathbb{E}(e_i | \mathcal{F}_0) = 0$ if $i > 0$, so

$$\begin{aligned}
\sum_{i=1}^{\infty} \mathbb{E}(f \circ T^i | \mathcal{F}_0) &= \sum_{k=1}^{\infty} \theta_k \sum_{i=1}^{N_k} e_{-N_k+i} \mathbb{1}_{T^{-i} A_k} \\
&= \sum_{k=1}^{\infty} \theta_k \sum_{i=1}^{N_k} e_{-N_k+i} \mathbb{1}_{A_k} + \sum_{k=1}^{\infty} \theta_k \sum_{i=1}^{N_k} e_{-N_k+i} (\mathbb{1}_{T^{-i} A_k \setminus A_k} - \mathbb{1}_{A_k \setminus T^{-i} A_k}).
\end{aligned} \tag{4}$$

Note that $|\mathbb{1}_{T^{-i} A_k \setminus A_k} - \mathbb{1}_{A_k \setminus T^{-i} A_k}| = \mathbb{1}_{T^{-i} A_k \Delta A_k}$ and by construction, $\mu(T^{-i} A_k \Delta A_k) \leq \varepsilon_k$ for $i \leq N_k$. Therefore,

$$\begin{aligned}
\mathbb{E}\left| \sum_{k=1}^{\infty} \theta_k \sum_{i=1}^{N_k} e_{-N_k+i} (\mathbb{1}_{T^{-i} A_k \setminus A_k} - \mathbb{1}_{A_k \setminus T^{-i} A_k}) \right| &\leq \sum_{k=1}^{\infty} \theta_k \sum_{i=1}^{N_k} \mathbb{E}|\mathbb{1}_{T^{-i} A_k \setminus A_k} - \mathbb{1}_{A_k \setminus T^{-i} A_k}| \\
&\leq \sum_{k=1}^{\infty} \theta_k N_k \varepsilon_k < \infty, \quad \text{by (3)}.
\end{aligned}$$

Hence, it remains to prove the $\mathbb{L}^1$-divergence of the first term in (4).

By disjointness of the $A_k$ and by independence between the $e_i$ and the $A_k$ for all $i$ and $k$,

$$\mathbb{E}\left| \sum_{k=1}^{\infty} \theta_k \sum_{i=1}^{N_k} e_{-N_k+i} \mathbb{1}_{A_k} \right| = \sum_{k=1}^{\infty} \theta_k \mathbb{E}\left| \sum_{i=1}^{N_k} e_{-N_k+i} \right| \mu(A_k).$$

Now, by independence of the $e_i$, we can use the Marcinkiewicz–Zygmund inequality (see e.g. Theorem 8.1 in [9] or Theorem 10.3.2 in [3]). There exists a constant $A > 0$, such that

$$\mathbb{E}\left|\sum_{i=1}^{N_k} e_{-N_k+i}\right| \geq A\mathbb{E}\left(\sum_{i=1}^{N_k} e_{-N_k+i}^2\right)^{1/2} = A\sqrt{N_k}.$$

Recall that $\mu(A_k) \geq a\rho_k$ (Lemma 2). So,

$$\mathbb{E}\left|\sum_{k=1}^{\infty} \theta_k \sum_{i=1}^{N_k} e_{-N_k+i}\mathbb{1}_{A_k}\right| \geq aA\sum_{k=1}^{\infty} \theta_k\sqrt{N_k}\rho_k = \infty, \quad \text{by assumption.}$$

This concludes the proof of the necessary condition.

**Sufficient condition.** We assume that $\sum_{k=1}^{\infty} \theta_k\sqrt{N_k}\rho_k < \infty$.

Let $I_n = \sum_{i=1}^{n} \mathbb{E}(f \circ T^i | \mathcal{F}_0)$. We will prove that $I_n \in \mathbb{L}^1$ for all $n$ and that the sequence $(I_n)_{n\geq 1}$ is a Cauchy sequence in $\mathbb{L}^1$. The proposition will follow from the completeness of $\mathbb{L}^1$.

To begin, we use the structure of the sets $A_k$ in the same way as in the first part of the proof. We have

$$
\begin{aligned}
I_n &= \sum_{i=1}^{n}\sum_{k=1}^{\infty} \mathbb{E}(f_k \circ T^i | \mathcal{F}_0)\mathbb{1}_{A_k} \circ T^i \\
&= \sum_{k=1}^{\infty} \theta_k \sum_{i=1}^{\min(n,N_k)} e_{-N_k+i}\mathbb{1}_{A_k} + \sum_{k=1}^{\infty} \theta_k \sum_{i=1}^{\min(n,N_k)} e_{-N_k+i}(\mathbb{1}_{T^{-i}A_k \setminus A_k} - \mathbb{1}_{A_k \setminus T^{-i}A_k}) \\
&= \Gamma_n^1 + \Gamma_n^2
\end{aligned}
\tag{5}
$$

and

$$
\begin{aligned}
\mathbb{E}|\Gamma_n^2| &\leq \sum_{k=1}^{\infty} \theta_k \sum_{i=1}^{\min(n,N_k)} \mathbb{E}|e_{-N_k+i}(\mathbb{1}_{T^{-i}A_k \setminus A_k} - \mathbb{1}_{A_k \setminus T^{-i}A_k})| \\
&\leq \sum_{k=1}^{\infty} \theta_k \sum_{i=1}^{\min(n,N_k)} \mu(T^{-i}A_k \Delta A_k) \\
&\leq \sum_{k=1}^{\infty} \theta_k N_k \varepsilon_k < \infty, \quad \text{by (3).}
\end{aligned}
\tag{6}
$$

On the other hand, by the Marcinkiewicz–Zygmund inequality, there exists a constant $B > 0$ such that

$$\mathbb{E}\left|\sum_{i=1}^{\min(n,N_k)} e_{-N_k+i}\right| \leq B\mathbb{E}\left(\sum_{i=1}^{\min(n,N_k)} e_{-N_k+i}^2\right)^{1/2}.$$

Recall that $\mu(A_k) \leq \rho_k$ (Lemma 2). Because $e_i$ is independent of $A_k$ for all $i$, $k$, we have

$$
\begin{aligned}
\mathbb{E}|\Gamma_n^1| &\leq B\sum_{k=1}^{\infty} \theta_k\mathbb{E}\left(\sum_{i=1}^{\min(n,N_k)} e_{-N_k+i}^2\right)^{1/2}\rho_k \\
&\leq B\sum_{k=1}^{\infty} \theta_k\sqrt{N_k}\rho_k < \infty, \quad \text{by assumption.}
\end{aligned}
\tag{7}
$$

Applying (6) and (7) to (5) shows that $I_n \in \mathbb{L}^1$ for all $n \in \mathbb{N}$.

Now, we will show that $(I_n)_{n \geq 1}$ is a Cauchy sequence in $\mathbb{L}^1$.

We fix $p \in \mathbb{N}^*$. We have

$$I_{n+p} - I_n = \sum_{k=1}^{\infty} \sum_{i=\min(n,N_k)+1}^{\min(n+p,N_k)} \mathbb{E}(f_k \circ T^i | \mathcal{F}_0) \mathbb{1}_{A_k} \circ T^i.$$

Using successively assumption (3), the Marcinkiewicz–Zygmund inequality and the independence between $e_i$ and $A_k$ for all $i$ and $k$ (see the calculus made before for $I_n$), we obtain

$$\mathbb{E}|I_{n+p} - I_n| \leq B \sum_{k=1}^{\infty} \theta_k \mathbb{E}\left( \sum_{i=\min(n,N_k)+1}^{\min(n+p,N_k)} e_{-N_k+i}^2 \right)^{1/2} \rho_k + \sum_{k=1}^{\infty} \theta_k \sum_{i=\min(n,N_k)+1}^{\min(n+p,N_k)} \mu(T^{-i} A_k \varDelta A_k).$$

Note that $\sum_{i=\min(n,N_k)+1}^{\min(n+p,N_k)}$ is empty if $N_k \leq n$, is composed of $N_k - (n+1)$ terms if $n < N_k < n+p$ and of $p$ terms otherwise. In the second and in the third case, the number of terms in the sum is less than $N_k$. So, for all $p \in \mathbb{N}^*$,

$$\mathbb{E}|I_{n+p} - I_n| \leq B \sum_{k:\, N_k > n} \theta_k \sqrt{N_k} \rho_k + \sum_{k:\, N_k > n} \theta_k N_k \varepsilon_k. \tag{8}$$

By assumption and hypothesis (3), $\sum_{k=1}^{\infty} \theta_k \sqrt{N_k} \rho_k < \infty$ and $\sum_{k=1}^{\infty} \theta_k N_k \varepsilon_k < \infty$. Hence, both sums in (8) go to 0 with $n \to \infty$ uniformly for all $p \in \mathbb{N}^*$.

$(I_n)_{n \geq 1}$ is thus a Cauchy sequence. $\qquad \square$

**Proposition 3.** *If $f \in \mathbb{L}^2$, the stationary process $(f \circ T^i)_{i \in \mathbb{Z}}$ verifies the projective criterion if and only if $\sum_{k=1}^{\infty} \theta_k^2 \sqrt{N_k} \rho_k < \infty$.*

**Proof.** It follows the idea of the proof of Proposition 2. So, some similar passages are given with less details.

**Necessary condition.** We assume that $\sum_{k=1}^{\infty} \theta_k^2 \sqrt{N_k} \rho_k = \infty$. We shall show that if $\sum_{i=1}^{\infty} f \mathbb{E}(f \circ T^i | \mathcal{F}_0)$ converges in $\mathbb{L}^1$, then $\mathbb{E}|\sum_{i=1}^{\infty} f \mathbb{E}(f \circ T^i | \mathcal{F}_0)| = \infty$. First,

$$\sum_{i=1}^{\infty} f \mathbb{E}(f \circ T^i | \mathcal{F}_0) = \sum_{k=1}^{\infty} f \theta_k \sum_{i=1}^{\infty} \mathbb{E}(e_{-N_k} \circ T^i | \mathcal{F}_0) \mathbb{1}_{T^{-i} A_k} = \sum_{k=1}^{\infty} f \theta_k \sum_{i=1}^{N_k} e_{-N_k+i} \mathbb{1}_{T^{-i} A_k}. \tag{9}$$

Like in the proof of Proposition 2, we decompose $\mathbb{1}_{T^{-i} A_k}$ into $\mathbb{1}_{A_k} + (\mathbb{1}_{T^{-i} A_k \setminus A_k} - \mathbb{1}_{A_k \setminus T^{-i} A_k})$. Applying the Cauchy–Schwarz inequality, we obtain

$$\mathbb{E}\left| \sum_{k=1}^{\infty} f \theta_k \sum_{i=1}^{N_k} e_{-N_k+i}(\mathbb{1}_{T^{-i} A_k \setminus A_k} - \mathbb{1}_{A_k \setminus T^{-i} A_k}) \right| \leq \|f\|_2 \sum_{k=1}^{\infty} \theta_k \sum_{i=1}^{N_k} \|\mathbb{1}_{T^{-i} A_k \varDelta A_k}\|_2$$

$$\leq \|f\|_2 \sum_{k=1}^{\infty} \theta_k N_k \sqrt{\varepsilon_k}.$$

Hypothesis (3), the fact that $f$ belongs to $\mathbb{L}^2$, and (9) show that the convergence of the integral

$$\mathbb{E}\left| \sum_{i=1}^{\infty} f \mathbb{E}(f \circ T^i | \mathcal{F}_0) \right|$$

is equivalent to the convergence of

$$\mathbb{E}\left|\sum_{k=1}^{\infty} f\theta_k\left(\sum_{i=1}^{N_k} e_{-N_k+i}\right)\mathbb{1}_{A_k}\right|.$$

Now, the sets $A_k$ being disjoint, we have

$$\sum_{k=1}^{\infty} f\theta_k\left(\sum_{i=1}^{N_k} e_{-N_k+i}\right)\mathbb{1}_{A_k} = \sum_{k=1}^{\infty}\sum_{j=1}^{\infty} \theta_j e_{-N_j}\theta_k\left(\sum_{i=1}^{N_k} e_{-N_k+i}\right)\mathbb{1}_{A_k}\mathbb{1}_{A_j}$$

$$= \sum_{k=1}^{\infty} \theta_k^2 e_{-N_k}\left(\sum_{i=1}^{N_k} e_{-N_k+i}\right)\mathbb{1}_{A_k}.$$

Using the disjointness of the sets $A_k$ again, the independence between $e_i$ and $A_k$ for all $i$ and $k$, the independence of the $e_{-N_k+i}$, $i = 0, \ldots, N_k$, the Marcinkiewicz–Zygmund inequality and the assumption, we obtain

$$\mathbb{E}\left|\sum_{k=1}^{\infty} f\theta_k\left(\sum_{i=1}^{N_k} e_{-N_k+i}\right)\mathbb{1}_{A_k}\right| = \sum_{k=1}^{\infty} \theta_k^2 \mathbb{E}|e_{-N_k}|\mathbb{E}\left|\sum_{i=1}^{N_k} e_{-N_k+i}\right|\mu(A_k)$$

$$\geq A\sum_{k=1}^{\infty} \theta_k^2 \mathbb{E}\left(\sum_{i=1}^{\min(n,N_k)} e_{-N_k+i}^2\right)^{1/2}\mu(A_k)$$

$$\geq aA\sum_{k=1}^{\infty} \theta_k^2\sqrt{N_k}\rho_k = \infty, \tag{10}$$

where $a$ comes from Lemma 2 and $A > 0$ comes from the Marcinkiewicz–Zygmund inequality.

**Sufficient condition.** We assume that $\sum_{k=1}^{\infty} \theta_k^2\sqrt{N_k}\rho_k < \infty$.

Let $J_n = \sum_{i=1}^{n} f\mathbb{E}(f \circ T^i|\mathcal{F}_0)$. We'll prove that $(J_n)_{n\geq 1}$ is a Cauchy sequence in $\mathbb{L}^1$, which proves the proposition.

First, we show that $J_n \in \mathbb{L}^1$ for all $n$, i.e., $\mathbb{E}|J_n| < \infty$ for all $n$. Indeed,

$$J_n = \sum_{k=1}^{\infty} f\theta_k\sum_{i=1}^{\min(n,N_k)} e_{-N_k+i}\mathbb{1}_{T^{-i}A_k}.$$

So, decomposing $\mathbb{1}_{T^{-i}A_k}$ into $\mathbb{1}_{A_k} + (\mathbb{1}_{T^{-i}A_k\setminus A_k} - \mathbb{1}_{A_k\setminus T^{-i}A_k})$, using the Cauchy–Schwarz inequality and (3), we show that it is enough to prove the convergence of

$$\mathbb{E}\left|\sum_{k=1}^{\infty} f\theta_k\left(\sum_{i=1}^{\min(n,N_k)} e_{-N_k+i}\right)\mathbb{1}_{A_k}\right|.$$

We repeat the calculus leading to (10) and we apply the Marcinkiewicz–Zygmund inequality. So, there exists $B > 0$ such that

$$\mathbb{E}\left|\sum_{k=1}^{\infty} f\theta_k\left(\sum_{i=1}^{\min(n,N_k)} e_{-N_k+i}\right)\mathbb{1}_{A_k}\right| = \sum_{k=1}^{\infty} \theta_k^2\mathbb{E}|e_{-N_k}|\mathbb{E}\left|\sum_{i=1}^{\min(n,N_k)} e_{-N_k+i}\right|\mu(A_k)$$

$$\leq B\sum_{k=1}^{\infty} \theta_k^2\sqrt{N_k}\rho_k < \infty, \quad \text{by assumption.}$$

Now, we fix $p \in \mathbb{N}^*$. By similar arguments, we can show that

$$\mathbb{E}|J_{n+p} - J_n| \leq B \sum_{k=1}^{\infty} \theta_k^2 \mathbb{E}\left(\sum_{i=\min(n,N_k)+1}^{\min(n+p,N_k)} e_{-N_k+i}^2\right)^{1/2} \rho_k + \|f\|_2 \sum_{k=1}^{\infty} \theta_k \sum_{i=\min(n,N_k)+1}^{\min(n+p,N_k)} \|\mathbb{1}_{T^{-i}A_k \triangle A_k}\|_2.$$

The same considerations about $\sum_{i=\min(n,N_k)+1}^{\min(n+p,N_k)}$ as in the proof of Proposition 2 lead to

$$\mathbb{E}|J_{n+p} - J_n| \leq B \sum_{k:\ N_k > n} \theta_k^2 \sqrt{N_k} \rho_k + \|f\|_2 \sum_{k:\ N_k > n} \theta_k N_k \sqrt{\varepsilon_k}.$$

By assumption and by (3), both sums go to 0 with $n \to \infty$, uniformly for $p \in \mathbb{N}^*$. Hence, $(J_n)_{n \geq 1}$ is a Cauchy sequence in $\mathbb{L}^1$, which is complete. $\qquad\square$

**Proposition 4.**

$$\sum_{n=1}^{\infty} \frac{\|\mathbb{E}(S_n(f)|\mathcal{F}_0)\|_2}{n^{3/2}} < \infty \quad \text{if and only if} \quad \sum_{n=1}^{\infty} \frac{(\sum_{k=1}^{\infty} \theta_k^2 \min(n,N_k)\rho_k)^{1/2}}{n^{3/2}} < \infty.$$

**Proof.** Note that

$$\mathbb{E}(S_n(f)|\mathcal{F}_0) = \sum_{k=1}^{\infty} \theta_k \sum_{i=1}^{\min(n,N_k)} e_{-N_k+i} \mathbb{1}_{A_k} + \sum_{k=1}^{\infty} \theta_k \sum_{i=1}^{\min(n,N_k)} e_{-N_k+i} (\mathbb{1}_{T^{-i}A_k \setminus A_k} - \mathbb{1}_{A_k \setminus T^{-i}A_k})$$

and

$$\left\| \sum_{k=1}^{\infty} \theta_k \sum_{i=1}^{\min(n,N_k)} e_{-N_k+i} (\mathbb{1}_{T^{-i}A_k \setminus A_k} - \mathbb{1}_{A_k \setminus T^{-i}A_k}) \right\|_2 \leq \sum_{k=1}^{\infty} \theta_k \min(n,N_k)\sqrt{\varepsilon_k}.$$

So, by (3),

$$\sum_{n=1}^{\infty} \frac{\|\mathbb{E}(S_n(f)|\mathcal{F}_0)\|_2}{n^{3/2}} < \infty \quad \text{if and only if} \quad \sum_{n=1}^{\infty} \frac{\|\sum_{k=1}^{\infty} \theta_k \sum_{i=1}^{\min(n,N_k)} e_{-N_k+i} \mathbb{1}_{A_k}\|_2}{n^{3/2}} < \infty.$$

Now, by independence, applying the Marcinkiewicz–Zygmund inequality, we can see that there exist $A$, $B > 0$ such that

$$aA \sum_{k=1}^{\infty} \theta_k^2 \min(n,N_k)\rho_k \leq \mathbb{E}\left| \sum_{k=1}^{\infty} \theta_k \sum_{i=1}^{\min(n,N_k)} e_{-N_k+i} \mathbb{1}_{A_k} \right|^2 \leq B \sum_{k=1}^{\infty} \theta_k^2 \min(n,N_k)\rho_k.$$

The proposition is proved. $\qquad\square$

## 5. Proof of the theorem

### 5.1. *Counterexample* 1, *proofs of* (i) *and* (iii)

In this section, we give an example of a function satisfying the projective criterion and also the Maxwell–Woodroofe condition but not the martingale-coboundary decomposition in $\mathbb{L}^1$. To do this, we consider the function $f$ defined at (2) by the sequences

$$\rho_k = \frac{1}{4^k}, \qquad N_k = 4^{2k} \quad \text{and} \quad \theta_k = \frac{1}{k}.$$

First,

$$\sum_{k=1}^{\infty} \theta_k^2 \rho_k = \sum_{k=1}^{\infty} \frac{1}{k^2 4^k} < \infty,$$

then, by Proposition 1, the function $f$ belongs to $\mathbb{L}^2$. We have

$$\sum_{k=1}^{\infty} \theta_k \sqrt{N_k} \rho_k = \sum_{k=1}^{\infty} \frac{1}{k} = \infty,$$

hence, by Proposition 2, the stationary process $(f \circ T^i)_{i \in \mathbb{Z}}$ does not admit a martingale-coboundary decomposition in $\mathbb{L}^1$. But,

$$\sum_{k=1}^{\infty} \theta_k^2 \sqrt{N_k} \rho_k = \sum_{k=1}^{\infty} \frac{1}{k^2} < \infty$$

and Proposition 3 show that it satisfies the projective criterion. This proves (i).

To verify that the process $(f \circ T^i)_{i \in \mathbb{Z}}$ satisfies the Maxwell–Woodroofe condition, by Proposition 4, we have to study the sums

$$\sum_{k=1}^{\infty} \theta_k^2 \min(n, N_k) \rho_k = \sum_{k=1}^{\lfloor (\ln n)/(2\ln 4) \rfloor} \theta_k^2 N_k \rho_k + n \sum_{k=\lfloor (\ln n)/(2\ln 4) \rfloor + 1}^{\infty} \theta_k^2 \rho_k. \tag{11}$$

The first term on the right-hand side can be estimated by

$$\sum_{k=1}^{\lfloor (\ln n)/(2\ln 4) \rfloor} \theta_k^2 N_k \rho_k = \sum_{k=1}^{\lfloor (\ln n)/(2\ln 4) \rfloor} \frac{4^k}{k^2} \leq \sum_{k=1}^{\lfloor (\ln n)/(2\ln 4) \rfloor} 4^k = \mathrm{O}(\sqrt{n}).$$

For the second term,

$$\sum_{k=\lfloor (\ln n)/(2\ln 4) \rfloor + 1}^{\infty} \theta_k^2 \rho_k = \sum_{k=\lfloor (\ln n)/(2\ln 4) \rfloor + 1}^{\infty} \frac{1}{k^2 4^k} \leq \sum_{k=\lfloor (\ln n)/(2\ln 4) \rfloor + 1}^{\infty} \frac{1}{4^k} = \mathrm{O}\left(\frac{1}{\sqrt{n}}\right).$$

From (11), we derive

$$\sum_{k=1}^{\infty} \theta_k^2 \min(n, N_k) \rho_k = \mathrm{O}(\sqrt{n})$$

and

$$\frac{(\sum_{k=1}^{\infty} \theta_k^2 \min(n, N_k) \rho_k)^{1/2}}{n^{3/2}} = \mathrm{O}(n^{-5/4}).$$

Therefore, by Proposition 4,

$$\sum_{n=1}^{\infty} \frac{\|\mathbb{E}(S_n(f)|\mathcal{F}_0)\|_2}{n^{3/2}} < \infty.$$

This proves (iii).

### 5.2. *Counterexample* 2, *proofs of* (ii) *and* (v)

Here, we show a process which satisfies the martingale-coboundary decomposition in $\mathbb{L}^1$ and the Maxwell–Woodroofe condition but fails to satisfy the projective criterion. We consider the function $f$ defined at (2), this time, by the sequences

$$\rho_k = \frac{1}{4^k}, \qquad N_k = k^2 \quad \text{and} \quad \theta_k = \frac{2^k}{k}.$$

We have

$$\sum_{k=1}^{\infty} \theta_k^2 \rho_k = \sum_{k=1}^{\infty} \frac{1}{k^2} < \infty, \qquad \sum_{k=1}^{\infty} \theta_k \sqrt{N_k} \rho_k = \sum_{k=1}^{\infty} \frac{1}{2^k} < \infty \quad \text{and} \quad \sum_{k=1}^{\infty} \theta_k^2 \sqrt{N_k} \rho_k = \sum_{k=1}^{\infty} \frac{1}{k} = \infty.$$

By Propositions 1, 2 and 3, $f$ belongs to $\mathbb{L}^2$ and satisfies the martingale-coboundary decomposition in $\mathbb{L}^1$ but not the projective criterion.

The process $(f \circ T^i)_{i \in \mathbb{Z}}$ verifies the Maxwell–Woodroofe condition. Indeed,

$$\sum_{k=1}^{\infty} \theta_k^2 \min(n, N_k) \rho_k = \sum_{k=1}^{\lfloor \sqrt{n} \rfloor} 1 + n \sum_{k=\lfloor \sqrt{n} \rfloor + 1}^{\infty} \frac{1}{k^2} = \mathrm{O}(\sqrt{n})$$

and, like in counterexample 1, using Proposition 4, we deduce

$$\sum_{n=1}^{\infty} \frac{\|\mathbb{E}(S_n(f)|\mathcal{F}_0)\|_2}{n^{3/2}} < \infty.$$

### 5.3. *Counterexample* 3, *proof of* (iv)

In this section our example verifies the martingale-coboundary decomposition in $\mathbb{L}^1$ with the invariance principle but not the Maxwell–Woodroofe condition. First, we consider the function $f$ defined at (2) by the sequences

$$\rho_k = \frac{1}{4^k}, \qquad N_k = 4^k \quad \text{and} \quad \theta_k = \frac{2^k}{k^{3/2}}.$$

We have

$$\sum_{k=1}^{\infty} \theta_k^2 \rho_k = \sum_{k=1}^{\infty} \frac{1}{k^3} < \infty \quad \text{and} \quad \sum_{k=1}^{\infty} \theta_k \sqrt{N_k} \rho_k = \sum_{k=1}^{\infty} \frac{1}{k^{3/2}} < \infty.$$

This implies that $f$ belongs to $\mathbb{L}^2$ and admits a martingale-coboundary decomposition in $\mathbb{L}^1$ (Propositions 1 and 2).

For the Maxwell–Woodroofe condition:

$$\sum_{k=1}^{\infty} \theta_k^2 \min(n, N_k) \rho_k = \sum_{k=1}^{\lfloor (\ln n)/\ln 4 \rfloor} \frac{4^k}{k^3} + n \sum_{k=\lfloor (\ln n)/\ln 4 \rfloor + 1}^{\infty} \frac{1}{k^3} \geq C \frac{n}{\ln^2 n},$$

for some $C > 0$. Therefore,

$$\sum_{n=1}^{\infty} \frac{(\sum_{k=1}^{\infty} \theta_k^2 \min(n, N_k) \rho_k)^{1/2}}{n^{3/2}} \geq \sqrt{C} \sum_{n=1}^{\infty} \frac{1}{n \ln n} = \infty$$

and by Proposition 4, the Maxwell–Woodroofe condition does not hold.

To prove (iv), it remains to show that the invariance principle holds. Actually, to do that, we will add hypotheses in the definition of the sets $A_k$. All the preceding results of this section will remain valid.

We have shown that $f$ admits a martingale-coboundary decomposition in $\mathbb{L}^1$. Thus, $f = m + g - g \circ T$, where $m, g \in \mathbb{L}^1$ and $(m \circ T^i)_{i \in \mathbb{Z}}$ is a martingale difference sequence. Here, we assume that $\mu(T^{-(N_k+1)}A_k \Delta A_k) \leq \varepsilon_k$ for all $k$ (in Lemma 2, take $N_k + 1$ instead of $N_k$). It is clear that this assumption does not change the previous results. Now, we can show:

**Proposition 5.** *In the decomposition $f = m + g - g \circ T$, $m$ belongs to $\mathbb{L}^2$.*

By the Billingsley and Ibragimov theorem for martingale difference sequences and by the stochastic boundedness of partial sums of $g - g \circ T$, it follows:

**Corollary 1.** *The process $(f \circ T^i)_{i \in \mathbb{Z}}$ verifies the CLT.*

**Proof.** Actually, we shall prove that $g - g \circ T \in \mathbb{L}^2$. In fact, see [20], $g = \sum_{i=0}^{\infty} \mathbb{E}(f \circ T^i | \mathcal{F}_0)$. So,

$$
\begin{aligned}
g - g \circ T &= \sum_{i=0}^{\infty} \mathbb{E}(f \circ T^i | \mathcal{F}_0) - \sum_{i=0}^{\infty} \mathbb{E}(f \circ T^{i+1} | \mathcal{F}_1) \\
&= \sum_{k=1}^{\infty} \theta_k \left( \sum_{i=0}^{N_k} e_{-N_k+i} \mathbb{1}_{T^{-i}A_k} - \sum_{i=0}^{N_k} e_{-N_k+i+1} \mathbb{1}_{T^{-(i+1)}A_k} \right) \\
&= \sum_{k=1}^{\infty} \theta_k (e_{-N_k} - e_1) \mathbb{1}_{A_k} + \sum_{k=1}^{\infty} \theta_k \sum_{i=0}^{N_k} e_{-N_k+i} (\mathbb{1}_{T^{-i}A_k \setminus A_k} - \mathbb{1}_{A_k \setminus T^{-i}A_k}) \\
&\quad - \sum_{k=1}^{\infty} \theta_k \sum_{i=0}^{N_k} e_{-N_k+i+1} (\mathbb{1}_{T^{-(i+1)}A_k \setminus A_k} - \mathbb{1}_{A_k \setminus T^{-(i+1)}A_k}).
\end{aligned}
$$

Now, by (3),

$$
\left\| \sum_{k=1}^{\infty} \theta_k \sum_{i=0}^{N_k} e_{-N_k+i} (\mathbb{1}_{T^{-i}A_k \setminus A_k} - \mathbb{1}_{A_k \setminus T^{-i}A_k}) \right\|_2 \leq \sum_{k=1}^{\infty} \theta_k \sum_{i=0}^{N_k} \| \mathbb{1}_{T^{-i}A_k \Delta A_k} \|_2 \leq \sum_{k=1}^{\infty} \theta_k N_k \sqrt{\varepsilon_k} < \infty. \tag{12}
$$

In the same way,

$$
\left\| \sum_{k=1}^{\infty} \theta_k \sum_{i=0}^{N_k} e_{-N_k+i+1} (\mathbb{1}_{T^{-(i+1)}A_k \setminus A_k} - \mathbb{1}_{A_k \setminus T^{-(i+1)}A_k}) \right\|_2 \leq \sum_{k=1}^{\infty} \theta_k N_k \sqrt{\varepsilon_k} < \infty. \tag{13}
$$

By disjointness of the sets $A_k$, by independence of the functions $e_{-N_k} - e_1$ and $\mathbb{1}_{A_k}$,

$$
\left\| \sum_{k=1}^{\infty} \theta_k (e_{-N_k} - e_1) \mathbb{1}_{A_k} \right\|_2^2 \leq \sum_{k=1}^{\infty} \theta_k^2 \| e_{-N_k} - e_1 \|_2^2 \mu(A_k) \leq 4 \sum_{k=1}^{\infty} \theta_k^2 \rho_k < \infty. \tag{14}
$$

(12)–(14) lead to the proposition. $\qquad \square$

Let $v_k = \theta_k \sum_{i=0}^{N_k} e_{-N_k} \circ T^i$. For $R \in \mathbb{N}^*$, the quantity $\sum_{k=1}^{R} \| v_k \|_\infty$ is finite. Thus, there exists an (not necessarily strictly) increasing sequence $(R_n)_{n \in \mathbb{N}} \to \infty$ such that

$$
\frac{1}{\sqrt{n}} \sum_{k=1}^{R_n - 1} \| v_k \|_\infty \xrightarrow[n \to \infty]{} 0. \tag{15}
$$

The sequence $(R_n)_{n \in \mathbb{N}}$ being fixed, we construct the sets $A_k$ in the following way. For all $k$, let $n_k$ be the greatest integer such that $R_{n_k} \leq k$. To define the sets $A_k$, we apply Lemma 2 with $(\max(n_k, N_k + 1))_{k \in \mathbb{N}}$ instead of $(N_k)_{k \in \mathbb{N}}$. Again, it is easy to see that previous results remain valid. With this construction, we have the following property.

$$\forall k \geq R_n, \forall i, j \in \{0, \ldots, n\}, \quad \mu(T^{-i} A_k \Delta T^{-j} A_k) \leq \varepsilon_k. \tag{16}$$

**Proposition 6.** *The process* $(f \circ T^i)_{i \in \mathbb{Z}}$ *verifies the invariance principle.*

**Proof.** Since $m \in \mathbb{L}^2$, as recalled in the Introduction, according to [21], it is enough to show that $\frac{1}{\sqrt{n}} \max_{i \leq n} |g \circ T^i| \underset{n \to \infty}{\longrightarrow} 0$ in probability. We have

$$g = \sum_{i=0}^{\infty} \mathbb{E}(f \circ T^i | \mathcal{F}_0) = \sum_{k=1}^{\infty} \theta_k \left( \sum_{i=0}^{N_k} e_{-N_k+i} \right) \mathbb{1}_{A_k} + \sum_{k=1}^{\infty} \theta_k \left( \sum_{i=0}^{N_k} e_{-N_k+i} (\mathbb{1}_{T^{-i}A_k \setminus A_k} - \mathbb{1}_{A_k \setminus T^{-i}A_k}) \right) = g_1 + g_2.$$

By the Markov inequality, for all $\lambda > 0$,

$$\mu \left\{ \max_{i \leq n} |g_2 \circ T^i| \geq \lambda \sqrt{n} \right\} \leq \frac{\mathbb{E}(\max_{i \leq n} |g_2 \circ T^i|)}{\lambda \sqrt{n}}$$

$$\leq \frac{1}{\lambda \sqrt{n}} \sum_{k=1}^{\infty} \theta_k (N_k + 1) \varepsilon_k \underset{n \to \infty}{\longrightarrow} 0, \quad \text{by (3)}.$$

So, $\frac{1}{\sqrt{n}} \max_{i \leq n} |g_2 \circ T^i|$ converges to 0 in probability.

It remains to prove the same thing for $g_1 = \sum_{k=1}^{\infty} v_k \mathbb{1}_{A_k}$. By (15),

$$\frac{1}{\sqrt{n}} \max_{i \leq n} \left| \sum_{k=1}^{R_n - 1} v_k \circ T^i \mathbb{1}_{A_k} \circ T^i \right| \leq \frac{1}{\sqrt{n}} \sum_{k=1}^{R_n - 1} \|v_k\|_{\infty} \underset{n \to \infty}{\longrightarrow} 0.$$

Hence, it converges to zero in probability.

Now, for all $\lambda > 0$, $\mu\{\max_{i \leq n} |\sum_{k=R_n}^{\infty} v_k \circ T^i \mathbb{1}_{A_k} \circ T^i| \geq 2\lambda \sqrt{n}\}$ is smaller than

$$\mu \left\{ \max_{i \leq n} \left| \sum_{k=R_n}^{\infty} v_k \circ T^i \mathbb{1}_{A_k} \right| \geq \lambda \sqrt{n} \right\} + \mu \left\{ \max_{i \leq n} \sum_{k=R_n}^{\infty} |v_k \circ T^i| \mathbb{1}_{T^{-i}A_k \Delta A_k} \geq \lambda \sqrt{n} \right\}.$$

For the first term, the Tchebychev inequality gives

$$\mu \left\{ \max_{i \leq n} \left| \sum_{k=R_n}^{\infty} v_k \circ T^i \mathbb{1}_{A_k} \right| \geq \lambda \sqrt{n} \right\} \leq \mu \left\{ \sum_{k=R_n}^{\infty} \max_{i \leq n} |v_k \circ T^i| \mathbb{1}_{A_k} \geq \lambda \sqrt{n} \right\}$$

$$\leq \frac{\mathbb{E}((\sum_{k=R_n}^{\infty} \max_{i \leq n} |v_k \circ T^i| \mathbb{1}_{A_k})^{1/3})}{\lambda n^{1/6}}$$

$$\leq \frac{\sum_{k=R_n}^{\infty} \theta_k^{1/3} (N_k + 1)^{1/3} \rho_k}{\lambda n^{1/6}} \underset{n \to \infty}{\longrightarrow} 0,$$

because

$$\sum_{k=1}^{\infty} \theta_k^{1/3} N_k^{1/3} \rho_k = \sum_{k=1}^{\infty} \frac{1}{\sqrt{k} 2^k} < \infty.$$

For the second term, the Markov inequality, assumptions (16) and (3) show convergence to zero with $n$.

Thus, $\frac{1}{\sqrt{n}} \max_{i \leq n} |g_1 \circ T^i| \longrightarrow_{n \to \infty} 0$ in probability. $\qquad\square$

(iv) is proved.

### 5.4. Counterexample 4, Proof of (vi)

To prove (vi), we have to improve our general model. To define our function $f$, in (2), we replace $f_k$ by $\theta_k \sum_{j=N_k+1}^{2N_k} e_{-j}$ and to define the sets $A_k$, we use Lemma 2 with $2N_k$ instead of $N_k$. Moreover, we assume that $\varepsilon_k$ is sufficiently small to have

$$\sum_{k=1}^{\infty} \theta_k N_k^2 \sqrt{\varepsilon_k} < \infty. \tag{17}$$

We shall see that if the sequences are well chosen, then the process $(f \circ T^i)_{i \in \mathbb{Z}}$ can verify the projective criterion but not the Maxwell–Woodroofe condition.

**Proposition 7.** *The function $f$ belongs to $\mathbb{L}^2$ if and only if $\sum_{k=1}^{\infty} \theta_k^2 N_k \rho_k < \infty$.*

**Proof.** It suffices to see that

$$\mathbb{E} \left| \sum_{k=1}^{\infty} \theta_k^2 \sum_{j=N_k+1}^{2N_k} e_{-j} \mathbb{1}_{A_k} \right|^2 = \sum_{k=1}^{\infty} \theta_k^2 \mathbb{E} \left| \sum_{j=N_k+1}^{2N_k} e_{-j} \right|^2 \mu(A_k)$$

and to use the independence between the $e_i$. $\qquad\square$

**Proposition 8.** *If $f \in \mathbb{L}^2$ and if $\sum_{k=1}^{\infty} \theta_k^2 N_k^2 \rho_k < \infty$, then $(f \circ T^i)_{i \in \mathbb{Z}}$ satisfies the projective criterion.*

**Proof.** Let $K_n = \sum_{i=1}^{n} f \mathbb{E}(f \circ T^i | \mathcal{F}_0)$. As in the proof of Proposition 3, using the properties of the sets $A_k$ and (17) we can see that $\|K_n\|_1 < \infty$ if and only if

$$\mathbb{E} \left| \sum_{k=1}^{\infty} \theta_k f \left( \sum_{i=1}^{n} \sum_{j=N_k+1}^{2N_k} \mathbb{E}(e_{i-j} | \mathcal{F}_0) \right) \mathbb{1}_{A_k} \right| < \infty.$$

Let us denote by $B_n^k$ the function $\sum_{i=1}^{n} \sum_{j=N_k+1}^{2N_k} \mathbb{E}(e_{i-j} | \mathcal{F}_0)$. Recall that $\mathbb{E}(e_i | \mathcal{F}_0)$ equals $e_i$ for $i \leq 0$ and equals 0 for $i > 0$, then

- for $n \geq 2N_k$,

$$B_n^k = N_k \sum_{j=0}^{N_k-1} e_{-j} + \sum_{j=0}^{N_k-1} (N_k - j) e_{-N_k-j};$$

- for $N_k < n < 2N_k$,

$$B_n^k = \sum_{j=0}^{2N_k-n} (n - N_k + j - 1) e_{-j} + \sum_{j=2N_k-n+1}^{N_k-1} N_k e_{-j} + \sum_{j=0}^{N_k-1} (N_k - j) e_{-N_k-j}; \text{ and}$$

- for $n \leq N_k$,

$$B_n^k = \sum_{j=0}^{n} (n - j) e_{-N_k+j} + \sum_{j=1}^{N_k} \min(n, N_k - j) e_{-N_k-j}.$$

In each case, by independence between the $e_i$, there exists $B > 0$ such that

$$\|B_n^k\|_2 \leq B N_k^{3/2}.$$

Thus, by the Cauchy–Schwarz inequality,

$$\mathbb{E}\left|\sum_{k=1}^{\infty} \theta_k f B_n^k \mathbb{1}_{A_k}\right| \leq \sum_{k=1}^{\infty} \theta_k^2 \left\|\sum_{j=N_k+1}^{2N_k} e_{-j}\right\|_2 \|B_n^k\|_2 \rho_k \leq B \sum_{k=1}^{\infty} \theta_k^2 N_k^2 \rho_k.$$

Therefore, $\sum_{k=1}^{\infty} \theta_k^2 N_k^2 \rho_k < \infty$ implies that $K_n$ belongs to $\mathbb{L}^1$ for all $n$. In the same manner, we can see that it also implies that $(K_n)_{n \in \mathbb{N}}$ is a Cauchy sequence in $\mathbb{L}^1$. The details are left to the reader. $\qquad \square$

Now, we choose the sequences in the definition of $f$. We take

$$\rho_k = \frac{1}{4^k}, \qquad N_k = 2^k \quad \text{and} \quad \theta_k = \frac{1}{k}.$$

Then

$$\sum_{k=1}^{\infty} \theta_k^2 N_k \rho_k = \sum_{k=1}^{\infty} \frac{1}{2^k k^2} < \infty \quad \text{and} \quad \sum_{k=1}^{\infty} \theta_k^2 N_k^2 \rho_k = \sum_{k=1}^{\infty} \frac{1}{k^2} < \infty.$$

Propositions 7 and 8 show that $f$ belongs to $\mathbb{L}^2$ and satisfies the projective criterion.

Using hypothesis (17) and the same observations as in the proof of Proposition 4, we see that the convergence of $\sum_{n=1}^{\infty} n^{-3/2} \|\mathbb{E}(S_n(f)|\mathcal{F}_0)\|_2$ is equivalent to the convergence of

$$\sum_{n=1}^{\infty} n^{-3/2} \left(\sum_{k=1}^{\infty} \theta_k^2 \mathbb{E}\left|\sum_{i=1}^{n} \sum_{j=N_k+1}^{2N_k} \mathbb{E}(e_{-j+i}|\mathcal{F}_0)\right|^2 \rho_k\right)^{1/2}.$$

For all $n \geq 2N_k$,

$$\sum_{i=1}^{n} \sum_{j=N_k+1}^{2N_k} \mathbb{E}(e_{-j+i}|\mathcal{F}_0) = N_k \sum_{j=0}^{N_k-1} e_{-j} + \sum_{j=0}^{N_k-1} (N_k - j) e_{-N_k-j},$$

and so

$$\mathbb{E}\left|\sum_{i=1}^{n} \sum_{j=N_k+1}^{2N_k} \mathbb{E}(e_{-j+i}|\mathcal{F}_0)\right|^2 \geq N_k^3.$$

Here, $N_k = 2^k$, so

$$\sum_{k=1}^{\infty} \theta_k^2 \mathbb{E}\left|\sum_{i=1}^{n} \sum_{j=N_k+1}^{2N_k} \mathbb{E}(e_{-j+i}|\mathcal{F}_0)\right|^2 \rho_k \geq \sum_{k=1}^{\lfloor (\ln n)/\ln 2 \rfloor - 1} \theta_k^2 N_k^3 \rho_k = \sum_{k=1}^{\lfloor (\ln n)/\ln 2 \rfloor - 1} \frac{2^k}{k^2} \geq C \frac{n}{\ln^2 n},$$

where $C$ is a positive constant. We derive that

$$\sum_{n=1}^{\infty} \frac{\|\mathbb{E}(S_n(f)|\mathcal{F}_0)\|_2}{n^{3/2}} = \infty,$$

i.e., the Maxwell–Woodroofe condition does not hold and (vi) is proved.

## Acknowledgment

The authors wish to express their thanks to the referee for several helpful comments and for pointing out a mistake in the proof of Proposition 2.